\documentclass{birkmult}

\usepackage{amsmath,amsthm}
\def\theequation{\thesection.\arabic{equation}}

\newtheorem{thm}{Theorem}[section]

\theoremstyle{definition}

\theoremstyle{remark}
\newtheorem{rem}[thm]{Remark}

\begin{document}

\title{Remark on spectral rigidity for magnetic Schr\"odinger
 operators.}

\author{Gregory Eskin} 
\address{  Department of Mathematics, UCLA,\\ Los Angeles,
CA 90095-1555, USA. }
\email{eskin@math.ucla.edu}
\author{James Ralston} 
\address{  Department of Mathematics, UCLA,\\ Los Angeles,
CA 90095-1555, USA. }
\email{ralston@math.ucla.edu}

\begin{abstract}
 We give a simple proof of Guillemin's theorem on the determination of
the magnetic field on the torus by the spectrum of the corresponding
Schr\"odinger operator.
\end{abstract}

\maketitle

\section{Introduction}

 This note is on inverse spectral theory for the Schr\"odinger
 operator on a flat two-dimensional torus with electric and magnetic
 potentials.
 This problem can be remarkably rigid. For generic flat tori, if the
 variation of the magnetic field is strictly less than its mean,
 and the total magnetic flux on the torus is $\pm 2\pi$, then the
 spectrum of the Schr\"odinger operator determines both the electric
 and magnetic fields.
 This is in marked contrast to both the Schr\"odinger operator
 without a magnetic field (see [3]) and the case of a magnetic
 field of mean zero (see [1]). In both those problems there are
 large families of isospectral fields, and rigidity results are
 much more difficult to obtain (see also [2]). The observation that
 there can be  spectral rigidity
 when the total flux is $\pm 2\pi$ is due to Guillemin ([5]). Here we
 give a short proof
 of the slightly stronger result stated above. Instead of
 thinking of
 the Hamiltonian as acting on functions with values in a line bundle
 over the torus $\Bbb R^2/L$, we
 think of the Hamiltonian as acting on functions
 on $\Bbb R^2$ which are invariant with respect to the \lq\lq
 magnetic translations" associated to $L$. However, these two settings
 are
 completely equivalent. Our assumption that the variation of the
 magnetic field $B(x)$ is strictly less than its mean $b_0$ takes the
 simple form  $|B(x)-b_0|<|b_0|$ for all $x$.

 The spectrum of the Laplacian plus lower order perturbations
 on flat tori has the feature that there are large families of spectral
 invariants corresponding to sets of geodesics with a fixed
 length. In analogy with results on $S^2$ Guillemin proposed the
 name \lq\lq band invariants" for these families. The nice feature of
 the problem
 discussed here is that only the simplest of the band invariants
 are needed to prove rigidity.

 The first complete solution of an inverse spectral problem
was Mark Krein's definitive analysis of the \lq\lq weighted
string", [9], [10]. Since that time many other inverse spectral
problems in one space dimension have been solved (see [11]). In
higher dimensions it is widely believed that, modulo natural
symmetries and deformations like gauge transformation,  most
problems will be spectrally rigid. However, so far there have been
relatively few settings where this has been proven (for instance
those in [6] and [14]) and many interesting examples where it
fails (see [12] and [4]). This should remain an active field of
research for many years to come, and one can reasonably say that
it began with the work of Mark Grigor'evich Krein. \vskip.2in

\section{Proof of Guillemin's Theorem}

\vskip.2in We begin with the smooth magnetic field $B$, periodic
with respect to the lattice $L$ in two dimensions, expanded in a
Fourier series in terms of the dual lattice $L^*$
$$B(x)=\sum_{\beta\in L^*}b_\beta e^{2\pi i\beta\cdot x}.$$
For this magnetic field we introduce the magnetic potential
$A=(A_1,A_2)$ with $\partial_{x_2}A_1-\partial_{x_1}A_2=B$, chosen
to be as periodic as possible, i.e.
$$A=A^0+A^1=\frac{b_0}{2}(x_2,-x_1)+\sum_{\beta\in L^*\backslash
 0}b_\beta
e^{2\pi i\beta\cdot x}{(\beta_2,-\beta_1)}{2\pi
i(\beta_1^2+\beta_2^2)}.$$ We also have a mean zero periodic
electric field which is the gradient of the mean zero periodic
potential
$$V(x)=\sum_{\beta\in L^*\backslash 0}v_\beta
e^{2\pi i\beta\cdot x}.$$ The quantum Hamiltonian for an electron
in these fields (with all physical constants set to 1)  is
$$H=(i\partial_x+A)^2+V.$$ \vskip.1in Let $D$ be a fundamental
domain for $L$. To define the domain of $H$ as an operator in
$L^2(D)$ we will use \lq\lq magnetic translation operators" (see
[13]). Letting $\{e_1,e_2\}$ and $\{e_1^*,e_2^*\}$ be a basis for
$L$ and the corresponding dual basis for $L^*$, define for
linearly independent vectors $v_1$ and $v_2$
$$T_ju(x)=e^{iv_j\cdot x}u(x+e_j),\qquad j=1,2.$$
Then the commutator $[T_1,T_2]$ is given by
$$[T_1,T_2]u(x)=(e^{iv_2\cdot e_1}-e^{iv_1\cdot
e_2})e^{i(v_1+v_2)\cdot x}u(x+e_1+e_2),$$ and the periodicity of
$A^1$ and $V$ implies that the commutator $[H,T_j]$ is given by
$$[H,T_j]u(x)=e^{iv_j\cdot
x}((i\partial_x+A(x)+A^0(e_j))^2-(i\partial_x+A(x)-v_j)^2)u(x+e_j).$$
Thus, in order for the $T_j$'s to commute with $H$ we require
$v_j=-A^0(e_j)$, and in order for the $T_j$'s to commute with each
other we require $A^0(e_1)\cdot e_2=-A^0(e_2)\cdot e_1= \pi l$ for
some integer $l$.  Note that this implies $A^0(e_1)=\pi l e^*_2$
and $A^0(e_2)=-\pi l e^*_1$ and $2\pi|l|=|b_0|\hbox{Area }(D)$,
and $b_0\hbox{Area }(D)=\int_DB(x)dx$ is the total magnetic flux.
Hence the assumption $b_0\neq 0$ is equivalent to nonzero flux,
and it implies $l\neq 0$. Defining the domain of $H$ to be the
subspace of $H^2(R^2)$ such that $T_j u=u,\ j=1,2$, we make $H$ a
self-adjoint operator in $L^2(D)$. \vskip.1in As in many previous
works we will look for spectral invariants for $H$ by studying the
wave trace. Letting $E(x,y,t)$ be the distribution kernel for the
fundamental solution for the initial value problem
$$u_{tt}+Hu=0\hbox{ in }\Bbb R^2_x\times \Bbb R_t, u(x,0)=f(x),\
u_t(x,0)=0,$$ the distribution kernel for the corresponding
initial value problem in $D\times \Bbb R_t$ is
\begin{equation}                   \label{eq:2.1}
E_D(x,y,t)=
\sum_{(m,n)\in \Bbb Z^2}T^m_1T^n_2E(x,y,t),
\end{equation} 
where the
operators $T_j$ act on the $x$ variable. Note that, since the
principal part of $\partial_t^2+H$ is $\partial_t^2-\Delta$,
$E(x,y,t)=0$ when $|x-y|^2>t^2$ and the sum in \eqref{eq:2.1} is has
 only a
finite number of nonzero terms for $t$ in a bounded interval. Thus
$[T_1,T_2]=0$ implies $T_jE_D(x,y,t)=E_D(x,y,t),\ j=1,2$. The
fundamental spectral invariant for this problem is the
distribution trace of the operator $E_D(t)$ corresponding to the
kernel $E_D(x,y,t)$. Conventionally (with all terms to be
interpreted in distribution sense) this is written
$$Tr(t)=\int_D E_D(x,x,t)dx.$$
\vskip.1in To avoid degeneracies in the contributions to $Tr(t)$
from the terms in (1), we assume that vectors in $L$ have distinct
lengths, i.e.
$$ d,d^\prime \in L \hbox{ and }|d|=|d^\prime| \hbox{ implies } d=\pm
 d^\prime.$$
Since $E(x,y,t)$ is singular as a distribution in $(x,y)$ only
when $|x-y|^2=t^2$, it now follows that the singularity of $Tr(t)$
at $t=|me_1+ne_2|$ comes from just two terms
\begin{equation}                              \label{eq:2.2}
\int_D[T^m_1T^n_2E(x,x,t)+T^{-m}_1T^{-n}_2E(x,x,t)]dx.
\end{equation}
To determine the spectral invariants coming from the leading terms
in the expansion of this singularity it is convenient to use the
Hadamard-H\"ormander expansion [7], [8] for $E(x,y,t)$. Beginning with
 the
forward fundamental solution, $E_+$, defined by
$(\partial_t^2+H)E_+=\delta(t)\delta(x-y)$ and $E_+=0$ for $t<0$
one has
\begin{equation}                   \label{eq:2.3}
E_+(t,x,y)\sim \sum_{\nu=0}^\infty
 a_\nu(x,y)e_\nu(t,|x-y|)
\end{equation}
where $e_\nu$ is chosen so that $(\partial_t^2-\Delta)e_\nu=\nu
e_{\nu-1}$ for $\nu>0$ and $e_0(t,|x-y|)$ is the forward
fundamental solution for $\partial_t^2-\Delta$. In two space
dimensions this means
 $$ e_\nu(t,|x-y|)=2^{-2\nu-1}\pi^{-1/2}{\mathcal
 X}^{\nu-1/2}_+(t^2-|x-y|^2)$$
 for $t>0$, $e_\nu=0$ for $t<0$. For $a>-1$ the distribution
 ${\mathcal X}^a_+$ is defined by ${\mathcal
 X}^a_+(s)=(\Gamma(a+1))^{-1}s^a$ for $s>0$ and ${\mathcal
 X}^a_+(s)=0$ for $s<0$. Hence the coefficients $a_\nu$
 are determined by the recursion
 $$
\nu a_nu+(x-y)\cdot \partial_xa_\nu
 -iA(x)\cdot(x-y)a_\nu+Ha_{\nu-1}=0,
 $$
 where $H$ acts in the variable $x$.
Solving this with the requirement that $a_0(y,y)=1$, we have

\begin{equation}                            \label{eq:2.4}
a_0(x,y)=\exp(i\int_0^1(x-y)\cdot A(y+s(x-y))ds)\hbox{ and
}
\end{equation}
\begin{equation}                             \label{eq:2.5} 
 a_1(x,y)=-a_0(x,y)(\int_0^1V(y+s(x-y))ds
+b(x,y),
\end{equation}
 where $b(x,y)$ is determined by $A(x)$.
 The fundamental solution $E(x,y,t)$ is given by
\begin{equation}                          \label{eq:2.6}
E(x,y,t)=\partial_t(E_+(t,x,y)-E_+(-t,x,y)).
\end{equation}
\vskip.1in We define
$$I(d)=\int_De^{-iA^0(d)\cdot x+i\int_0^1d\cdot A(x+sd)ds}dx
= \int_De^{i(2A^0(x)\cdot d-\int_0^1d\cdot A^1(x+sd)ds)}dx, $$ and
$$J(d)=\int_D[\int_0^1(V(x+sd)+b(x+sd,x))ds]e^{i(2A^0(x)\cdot
 d-\int_0^1d\cdot A^1(x+sd)ds)}dx.$$ From \eqref{eq:2.2}-\eqref{eq:2.6}
 one 
sees that $I(d)+I(-d)$ and
$J(d)+J(-d)$ are spectral invariants for $H$. However, the
periodicity implies that $I(d)=I(-d)$ and $J(d)=J(-d)$. \vskip.1in
The rest of this article is devoted to studying $I(d)$ and $J(d)$.
We have $d=me_1+ne_2=k(m_0e_1+n_0e_2)$, $k\in \Bbb N$ and
gcd($m_0,n_0$)=1. Let $\delta=-n_0e^*_1+m_0e^*_2$. Then we have
$$A^0(d)=\frac{b_0}{2}(d_2,-d_1)=\pi kl\delta.$$
Since $\int_0^1e^{2\pi is\beta\cdot d }ds=0$ when $\beta\cdot
d\neq 0$, the terms in the Fourier series for $A^1$ which
contribute to $I(d)$ have $\beta\cdot d=0$, and this implies
$$\beta =p\delta
=\frac{pb_0}{2\pi kl}(d_2,-d_1),\ p\in \Bbb Z\backslash 0.$$ Hence,
$d\cdot (\beta_2,-\beta_1)(2\pi
i(\beta_1^2+\beta_2^2))^{-1}=ikl(pb_0)^{-1}$, and $I(d)$ reduces
to
$$\int_D \exp(2\pi ikl(-\delta\cdot x+\sum_{p\in \Bbb Z\backslash
0}\frac{ib_{p\delta}}{2\pi pb_0}e^{2\pi ip\delta\cdot x}))dx.$$
Defining
$$B_\delta(s)=\sum_{p\in \Bbb Z\backslash
0}b_{p\delta}e^{2\pi ips}\hbox{ and } A^1_\delta(s)= \sum_{p\in
\Bbb Z\backslash 0}\frac{b_{p\delta}}{2\pi ip}e^{2\pi i ps}$$ (note
that $\frac{d}{ds}A^1_\delta(s)=B_\delta(s)$), we have
$$I(d)=\int_D\exp(-i2\pi kl(\delta\cdot x+\frac{1}
 {b_0}A^1_\delta(\delta\cdot x))dx.$$
 Extending $\delta$ to a basis for
 $L^*,\ \{\delta, \delta^\prime\}$, and letting
 $\{\gamma,\gamma^\prime\}$ be the dual basis for $L$, we make the
 change of variables $x=s\gamma +u\gamma^\prime$, and choose
 $$D=\{s\gamma +u\gamma^\prime: 0\leq s,u\leq 1\}.$$ Then
 we have
 $$I(d)=c(d)\int_0^1\exp(-2\pi ikl(s+\frac{1}
 {b_0}A^1_\delta(s))ds,$$
 where $c(d)$ is the Jacobian factor, and only depends on $d$.
 Since we have this
 spectral invariant for all $k\neq 0$, it
 follows that
\begin{equation}                         \label{eq:2.7}
 \int_0^1f(s +\frac{1}{b_0}A^1_\delta(s))ds 
\end{equation}
 is a spectral invariant for any function $f$ which can be expanded
 in terms of $\{e^{-2\pi ikly}\}_{k\in \Bbb Z}$, i.e. for any $f\in
 L^2_{loc}(\Bbb R)$
 which has period $1/l$.

\begin{thm}                           
\label{theo:2.1}
Assume that $l=1$ and $|b_0|>\max |B(x)-b_0|$.
Then the spectrum of $H$ determines $B$.
\end{thm}

\begin{rem}
Since $b_0$ is the average of $B(x)$ on a fundamental
domain, the hypothesis here is a constraint on how much $B$ varies
instead of constraint on the size of $B$. 
\end{rem}
\begin{rem} Since we assume $1=l=A^0(e_1)\cdot e_2/\pi
=|b_0|\hbox{Area}(D)/2\pi$ for a fundamental domain $D$, this
assumption fixes $|b_0|$ when $L$ is fixed.
\end{rem}

\begin{proof}[Proof of Theorem~\ref{theo:2.1}] Since
$B_\delta(x)=\int_0^1(B(x+sd)-b_0)ds$, the hypotheses imply that
the derivative of $s+A_\delta^1(s)/b_0$ is strictly positive and
the inverse function $s(y)$ to $y=s+A^1_\delta(s)/b_0$ is defined
on the range of $s+A^1_\delta(s)/b_0$ for $s\in [0,1]$. Since
$A^1_\delta$ has period 1, the range is
$I=[A^1_\delta(0)/b_0,A^1_\delta(0)/b_0+1]$. Letting $f$ in
 \eqref{eq:2.7}
tend to the $\delta$-function at $y$, the limit of \eqref{eq:2.7} is
$s^\prime(y)$ if $y=y(s)$ for $s\in [0,1]$. If $y\neq y(s)$ for
$s\in [0,1]$ then the limit of \eqref{eq:2.7} is $s^\prime(y^*)$, where
$y^*\in I$, and $y^*=y \hbox{ mod }1$.  In other words taking
these limits we recover a function of period 1 in $y$ which agrees
with $s^\prime(y)$ on $I$. Thus we recover $A^1_\delta(s)$ modulo
an additive constant, and we obtain $B_\delta(s)$ by taking the
derivative. Since we can carry out this argument for all prime
elements $\delta\in L^*$, we recover the full Fourier expansion of
$B$. 
\end{proof}

We now turn to the recovery of $V$. The preceding
analysis shows that, keeping the same $d\in L$ as above, the
spectral invariant $J(d)$, modulo terms determined by $A(x)$,
reduces to
\begin{equation}                       \label{eq:2.8}
\tilde J(d)=c(d)\int_0^1V_\delta(s)e^{-2\pi
 ik(s+A_\delta(s))}ds,
\end{equation}
where
$$V_\delta(s)=\sum_{p\in \Bbb Z\backslash
0}v_{p\delta}e^{2\pi ips}.$$
This immediately gives the following:

\begin{thm} \label{theo:2.4}
Under the hypotheses of Theorem~\ref{theo:2.1} the spectrum of
$H$ determines $V$. 
\end{thm}
\begin{proof}[Proof of Theorem~\ref{theo:2.4}.]
Since we are assuming the hypotheses of Theorem~\ref{theo:2.1}, we have
 the
function $s(y)$ and can make the substitution $s=s(y)$ in
 \eqref{eq:2.8}.
That gives
$$\tilde
J(d)=c(d)\int_{A^1_\delta(0)}^{A^1_\delta(0)+1}V_\delta(s(y)e^{-2\pi
iky}s^\prime(y)dy,$$ but, since $y(s+1)=y(s)$, we can extend
$s(y)$ smoothly to the whole line by defining $s(y+1)=s(y)+1$.
Thus, since $V_\delta(s)$ has period 1 in $s$, we have
$$\tilde J(d)=c(d)\int_0^1V_\delta(s(y))s^\prime(y)e^{-2\pi
iky}dy.$$ Since we have this spectral invariant for $k\in {\Bbb
Z}\backslash 0$, we recover the Fourier series of
$V_\delta(s(y))s^\prime(y)$, and, hence, since $s(y)$ is
determined by $A_\delta^1(s)$, we have $V_\delta(s)$. As before,
since we can carry out this argument for all prime elements
$\delta\in L^*$, we recover the full Fourier expansion of $V$.
\end{proof}

\begin{rem} If $l=p/q,\ p,q\in \Bbb N$, for the
lattice $L$, then $l=1$ for the lattice $L_0$ generated by
$c_1e_1+c_2e_2$ and $d_1e_1+d_2$ when $p(c_1d_2-c_2d_1)=q$. So if
$B(x)$ and $V(x)$ are periodic with respect to $L_0$, Theorems 1
and 2 apply in the sense that the spectrum of $H$ on the torus
$\Bbb R^2/L_0$ determines $B(x)$ and $V(x)$. Note that $B(x)$ and
$V(x)$ will automatically be periodic with respect to $L_0$ when
$l=1/q$.
\end{rem}

\end{document}